\documentstyle[12pt]{article}
\newcommand{\too}{\longrightarrow}
\newcommand{\om}{\omega}
\newcommand{\Om}{\Omega}
\newcommand{\na}{\nabla}
\newcommand{\wi}{\widetilde}
\newcommand{\al}{\alpha}
\newcommand{\be}{\beta}

\newcommand{\Ga}{\Gamma}

\newcommand{\entraine}{\Longrightarrow}
\newcommand{\inj}{\hookrightarrow}

\newcommand{\ssi}{\Longleftrightarrow}
\def \reel{ {\rm I}\!{\rm R} }

 \def \rat{ {\rm Q}\kern-.65em {}^{{}_/ }}

 \def\ent{{{\rm Z}\mkern-5.5mu{\rm Z}}}

\newtheorem{th}{Theorem}[section]
\newtheorem{pr}{Proposition}[section]

\newtheorem{co}{Corollary}[section]
\title{ The modular class of a regular Poisson manifold
and the Reeb invariant of its symplectic foliation}
\author{Abdelhak Abouqateb and Mohamed Boucetta} \date{}\parindent=0cm \begin{document} \maketitle

{\bf Abstract.} \footnote[1]{The second authors' research was
supported by The Third World Academy of Sciences RGA No 01-301
RG/Maths/AC}

We show that, for any regular Poisson manifold, there is an
injective  natural linear map from the  first leafwise cohomology
space into the first Poisson cohomology space which maps the Reeb
class of the symplectic foliation to the modular class of the
Poisson manifold. The Riemannian interpretation of those classes
will permit us to show that a regular Poisson manifold whose
symplectic foliation is of codimension one is unimodular if and
only if its symplectic foliation is Riemannian foliation. It
permit us also to construct examples of unimodular Poisson
manifolds and other which are not unimodular. Finally, we prove
that the first leafwise cohomology space is an invariant of Morita
equivalence.

{\it Key words.}   Poisson manifold, Modular class, Riemannian
foliation.

{\it 2000 M. S.  C. Primary 53D17; Secondary 53C12.}
\section{Introduction}
The Reeb invariant of a foliated manifold is an obstruction lying
in the first leafwise cohomology to the existence of a volume
normal form invariant by the vector fields tangent to the
foliation [5]. The modular class of a Poisson manifold was
introduced by Weinstein [11]. It is an obstruction lying in the
first Poisson cohomology to the existence of a volume form
invariant with respect to the Hamiltonian flows.  For a regular
Poisson manifold, Weinstein [11 ] pointed out that the two classes
are closely related without giving  an explicit relation between
them. In fact, the two classes represent the same mathematical
objet. We will show that the first leafwise cohomology space is,
in natural way, a subspace of the first Poisson cohomology space
and the Reeb class of the symplectic foliation agrees  with the
modular class of the Poisson manifold.

On a Riemannian foliated manifold,  we remark that  the tangent
mean curvature gives arise to a tangential 1-forme whose leafwise
cohomology class is the Reeb invariant. This remark and the fact
that the Reeb invariant is the same object as the modular class
permit us to have the following corollaries.
\begin{co}Let $(P,\pi)$ be a regular Poisson manifold. If the
symplectic foliation is Riemannian then the modular class of $P$
vanishes.\end{co}
\begin{co}Let $(P,\pi)$ be a regular Poisson manifold for which the symplectic foliation is
 transversally oriented of codimension 1. The following assertions are equivalent:

1) The modular class of $P$ vanishes.

2) The symplectic foliation is Riemannian.\end{co}
\begin{co} Let $(P,\pi)$ a simply connected and compact regular Poisson manifold for which
the symplectic foliation is transversally oriented of codimension
1. Then $mod(P)\not=0$.\end{co}

With this corollary in mind, we construct many examples of regular
Poisson manifold with vanishing modular class and many examples
with non-vanishing modular class.

The first Poisson cohomology spaces of Morita equivalent Poisson
manifolds are isomorphic, according to Ginzburg an Lu [4], and
Ginzburg [3] has shown that the modular classes are compatible
with this isomorphism.  We will show that the first leafwise
cohomology spaces are also compatible with this isomorphism.

\section{The Reeb class of a foliation and its Riemannian
interpretation}

Let $M$ be a differentiable manifold endowed with a transversally
oriented foliation $\cal F$ of dimension $p$ and of codimension
$q$. We denote $T{\cal F}$ the tangent bundle to the foliation and
${\cal X}({\cal F})$ the space of vector fields tangent to $\cal
F$. Let ${\cal A}^r_{\cal F}$ denote the space of sections of the
bundle $\wedge^rT^*{\cal F}\too M$. The elements of ${\cal
A}^r_{\cal F}$ are called tangential differential $r$-forms. The
expression
\begin{eqnarray*}
d_{\cal F}\al(X_1,\ldots,X_{r+1})&=&\sum_{i=1}^{r+1}(-1)^{i+1}X_i.\al(X_1,\ldots,\hat X_i,\ldots,X_{r+1})\\
&+&\sum_{i<j}(-1)^{i+j}\al([X_i,X_j],X_1,\ldots,\hat
X_i,\ldots,\hat X_j,\ldots,X_{r+1}),\qquad(1)\end{eqnarray*}where
$\al\in{\cal A}^r_{\cal F}$ and $X_1,\ldots,X_{r+1}\in{\cal
X}({\cal F})$, defines a degree one differential operator $d_{\cal
F}$ that satisfies $d_{\cal F}^2=0$. The induced cohomology
$H_{\cal F}^*(M)$ is called the leafwise cohomology.

An orientation of the normal bundle to $\cal F$ is a differential
$q$-form $\nu$ on $M$ such that $\nu_x\not=0$ for any $x\in M$ and
such that $i_X\nu=0$ for any $X\in{\cal X}({\cal F})$.

For any $X\in{\cal X}({\cal F})$, $L_X\nu$ is proportional to
$\nu$ and one can define a tangential 1-form $\al_{\cal F}$ by
$$L_X\nu=\al_{\cal F}(X)\nu.\eqno(2)$$

From $L_{[X,Y]}\nu=L_X\circ L_Y\nu-L_Y\circ L_X\nu$, we have
$$d_{\cal F}\al_{\cal F}=0.\eqno(3)$$
The  cohomology class of $\al_{\cal F}$  denoted by $mod({\cal
F})$ is the Reeb class of the foliation.

The normal bundle to $\cal F$ carries an orientation $\nu$ such
that $L_X\nu=0$ for any $X\in{\cal X}({\cal F})$ if and only if
$mod({\cal F})=0$ (see [5]).

Now, we give the Riemannian interpretation of the Reeb class.

Let $g$ be a Riemann metric on $M$ and let $\na$ be the associated
Levi-Civita connection. We denote $T^\perp{\cal F}$ the
orthogonal distribution to  $T\cal F$ and
 ${\cal X}({\cal F}^\perp)$ the space of vector fields
tangent to $T^\perp{\cal F}$ . For any vector field $X$, denote
$X^{\cal F}$ its component in ${\cal X}({\cal F})$ and $X^{{\cal
F}^\perp}$ its component in ${\cal X}({\cal F}^\perp)$. The
orthogonal volume form to the foliation is  the differential
$q$-form $\eta$ defined by $\eta(Y_1,\ldots,Y_q)=1$ for any
orthonormal oriented frame  $(Y_1,\ldots,Y_q)$ in ${\cal X}({\cal
  F}^\perp)$ and $i_X\eta=0$ for any $X\in{\cal X}({\cal F})$.

A straightforward calculation gives
$$L_X\eta(Y_1,\ldots,Y_q)=\sum_{i=1}^qg(\na_{Y_i}X,Y_i)=-\sum_{i=1}^qg(\na_{Y_i}Y_i,X).\eqno(4)$$
The  second  fundamental form of   $T^\perp{\cal F}$   is the
tensor field $B^\perp: {\cal X}({\cal F}^\perp)\times{\cal
X}({\cal F}^\perp)\too{\cal
  X}({\cal F})$ given by
$$B^\perp(Y_1,Y_2)=\frac12[\na_{Y_1}Y_2+\na_{Y_2}Y_1]^{\cal
F}.$$Its  trace with respect to  $g$, called the tangent mean
curvature, is a vector field $H^\perp$ tangent to $\cal F$. We
define a tangential 1-form  $K^\perp\in{\cal
  A}^1_{\cal F}$ by
$$K^\perp(X)=g(X,H^\perp),\qquad X\in {\cal X}({\cal F}).\eqno(5)$$

(3) can be written
$$L_X\eta(Y_1,\ldots,Y_q)=-K^\perp(X)\eqno(6)$$and so
$$mod({\cal F})=-[K^\perp].\eqno(7)$$

{\bf Remark.} This formula can be compared to the metric formula for the Godbillon-Vey invariant ( see [8]).

\begin{pr}Let $M$ be a  differentiable manifold endowed with a transversally
oriented foliation $\cal F$ of dimension $p$ and of codimension
$q$. The following assertions are equivalent:

1) The normal bundle of $\cal F$ carries an orientation invariant
by the vector fields tangent to the foliation.

2) $mod({\cal F})=0.$

3) There is a Riemann metric on $M$ with vanishing tangent mean
curvature.\end{pr}

{\bf Proof:} We have shown  that $1)\ssi 2)$and
   that $3)\entraine 2)$. We show now  $2)\entraine 3)$. Let $g$ be a
   Riemann metric on $M$.
If $mod({\cal F})=0$
   then
   there is a smooth function  $h\in C^\infty(M)$ such that $d_{\cal
   F}h=K^\perp$. It's easy to verify that the tangent mean curvature of the Riemann metric
   $g_1=e^{-\frac{2h}{q}}g$
   vanishes.$\Box$

It's known that the   second  fundamental form  $B$ vanishes if
and only if   $g$ is bundle-like [9] so we get the following
proposition.
\begin{pr} Let $M$ be a  differentiable manifold endowed with a transversally
oriented foliation $\cal F$  of codimension $1$. The following
assertions are equivalent:

1) $mod({\cal F)}=0$.

2) $\cal F$ is a Riemannian foliation.

3) $\cal F$ is defined by a closed 1-form.\end{pr}
 \section{The modular class of a Poisson manifold}
Many fundamental definitions and results about Poisson manifolds
can be found in Vaisman's monograph [10].

Let $P$ be a Poisson manifold with Poisson tensor $\pi$. We  have
a bundle map $\pi:T^*P\too TP$ defined by
$$\be(\pi(\al))=\pi(\al,\be),\qquad\al,\be\in T^*P.\eqno(8)$$
On the space of differential 1-forms $\Om^1(P)$, the Poisson
tensor induces a Lie bracket
\begin{eqnarray*}
[\al,\be]_{\pi}&=&L_{\pi(\al)}\be-L_{\pi(\be)}\al-d(\pi(\al,\be))\\
&=&i_{\pi(\al)}d\be-i_{\pi(\be)}d\al+d(\pi(\al,\be)).\qquad\qquad\qquad\qquad\qquad\quad
(9) \end{eqnarray*}

For this Lie bracket and the usual Lie bracket on vector fields,
the bundle map $\pi$ induces a Lie algebra homomorphism
$\pi:\Om^1(P)\too{\cal X}(P)$:
$$\pi([\al,\be]_{\pi})=[\pi(\al),\pi(\be)].\eqno(10)$$

The Poisson cohomology of a Poisson manifold $(P,\pi)$ is the
cohomology of the chain complex $({\cal X}^*(P),d_\pi)$ where, for
$0\leq p\leq dimP$, ${\cal X}^p(P)$ is the
$C^\infty(P,\reel)$-module of $p$-multi-vector fields and $d_\pi$
is given by
\begin{eqnarray*}
d_\pi Q(\al_0,\ldots,\al_p)&=&\sum_{j=0}^p(-1)^j
\pi(\al_j).Q(\al_0,\ldots,\hat\al_j,\ldots,\al_p)\\
&+&\sum_{i<j}(-1)^{i+j}
Q([\al_i,\al_j]_\pi,\al_0,\ldots,\hat\al_i,\ldots,\hat\al_j,\ldots,\al_p).\qquad(11)
\end{eqnarray*}
We denote $H_\pi^*(P)$ the spaces of cohomology.

The modular class of $(P,\pi)$ is the obstruction to the existence
of a volume form on $P$ which is invariant with respect to
Hamiltonian flows. More explicitly, let $\mu$ be a volume form on
$P$. As shown in [11], the operator $\phi_\mu:f\mapsto
div_\mu\pi(df)$ is a derivation and hence a vector field called
the modular vector field of $(P,\pi)$ with respect to the volume
form $\mu$. Also $$L_{\phi_\mu}\pi=0\qquad\mbox{ and}\qquad
L_{\phi_\mu}\mu=0.\eqno(12)$$

If we replace $\mu$ by $a\mu$, where $a$ is a positive function,
the modular vector fields becomes
$$\phi_{a\mu}=\phi_\mu+\pi(d(Loga)).\eqno(13)$$Thus the first Poisson
cohomology class of $\phi_\mu$ is independent of $\mu$, we call it
the modular class of $(P,\pi)$ and we denote it $mod(P)$. The
Poisson manifold is unimodular if its modular class vanishes.
\section{ Link between the Reeb class and the modular class of a
regular Poisson manifold}

Let $(P,\pi)$ be a Poisson manifold whose symplectic foliation,
denoted by $\cal F$, is a regular foliation transversally oriented
of dimension $2p$ and of codimension $q$. As shown in[11,pp. 385],
the modular vector field of $P$ is closely related to the Reeb
tangential 1-form. More explicitly, let $\om\in\Ga(\wedge^
2T^*{\cal F})$ be the leafwise symplectic form given
by$$\om(u,v)=\pi(\pi^{-1}(u),\pi^{-1}(v)), \qquad u,v\in T{\cal
F}.$$ Let $\nu$ be a transverse orientation of $\cal F$ i.e a
differential $q$-form on $P$ such that $i_{\pi(df)}\nu=0$ for any
smooth function $f$ and $\nu_x\not=0$ for any $x\in P$. Choose $F$
a supplement distribution to $T\cal F$ and extend $\om$ to a
differential 2-form on $P$ by setting $i_X\om=0$ for any $X$
tangent to $F$. The form $\mu=\wedge^p\om\wedge\nu$  is a volume
form on $P$.

For any $f\in C^\infty(P)$,  since
$[L_{\pi(df)}(\wedge^p\om)]\wedge\nu=0$,
 we have
$$
L_{\pi(df)}\mu=(\wedge^p\om)\wedge L_{\pi(df)}\nu=\al_{\cal
F}(\pi(df))\mu=-\pi(\al)(f)\mu$$where $\al$ is any differential
1-form on $P$ whose restriction to $\cal F$ is $\al_{\cal F}$.

 $\pi(\al)$ depends only on $\al_{\cal F}$ and we will denote it by  $\pi(\al_{\cal F})$. We get
$$\phi_{(\wedge^p\om)\wedge\nu}=-\pi(\al_{\cal F}).\eqno(14)$$
In the following, we will give a more precise interpretation of
this relation.

For $1\leq p\leq dimP$, we consider the subspace ${\cal
X}^p_0(P)\subset{\cal X}^p(P)$ of $p$-multi-vector fields $Q$ such
that $i_\al Q=0${ for any} $\al\in Ker\pi.$ It's easy to verify
that $d_\pi({\cal X}^p_0)\subset{\cal X}^{p+1}_0$. The natural
injection ${\cal X}^p_0\inj{\cal X}^p$ induces a linear map
$H^*({\cal X}^p_0)\too H_\pi^*(P)$ which is  injective for $*=1$.

Let $\pi: {\cal A}_{\cal F}^p(P)\too{\cal X}^p_0(P)$ be the map
given by
$$\pi(\om)(\al_1,\ldots,\al_p)=\om(\pi(\al_1),\ldots,\pi(\al_p)).$$
It is easy to verify that $\pi$ is an isomorphism and $\pi(d_{\cal
F}\om)=d_\pi\pi(\om)$ and hence $\pi$ induces an isomorphism
$$\pi^*:H^p_{\cal F}(P)\too
H^p({\cal X}^p_0(P)).\eqno(15)$$ So, we have shown the following
theorem.
\begin{th} Let  $(P,\pi)$ be a regular Poisson manifold for which the symplectic foliation
is transversally oriented.  $\pi$ induces a linear  injection
$$\pi^*:H_{\cal F}^1(P)\inj H_\pi^1(P)\eqno(16)$$
 and we have
$$\pi^*(mod({\cal F}))=mod(P).\eqno(17)$$\end{th}
{\bf Remark.} The fact that the leafwise cohomology spaces embed
in the Poisson cohomology is known ( see [10]).

 The following corollaries  are a consequence of this theorem,
Proposition 2.1 and Proposition 2.2.
\begin{co}Let $(P,\pi)$ be a regular Poisson manifold. If the
symplectic foliation is Riemannian then the modular class of $P$
vanishes.\end{co}
\begin{co} Let  $(P,\pi)$ be a regular Poisson manifold for which the symplectic foliation
is transversally oriented. The following assertions are
equivalent:

1) $P$ is unimodular.

2) The symplectic foliation carries an invariant normal volume
form.

3) There is a Riemannian metric $g$ on $P$ such that
$\pi(K^\perp)$ is a hamiltonian vector fields.\end{co}

\begin{co} Let  $(P,\pi)$ be a regular Poisson manifold for which the symplectic foliation is
 transversally oriented of codimension 1. The following assertions are equivalent:

1) $P$ is unimodular.

2) The symplectic foliation is Riemannian.\end{co}
\begin{co} Let $(P,\pi)$ a simply connected and compact regular Poisson manifold for which
the symplectic foliation is transversally oriented of codimension
1. Then $mod(P)\not=0$.\end{co}
\section{Examples}
According to the above sections, we will give some illustration
examples of regular Poisson manifolds with vanishing modular class
( Reeb class)  and other with non-vanishing modular class ( Reeb
class).

1. Any  oriented foliation of codimension 1 on the sphere $S^3$ is
the symplectic foliation of a Poisson structure on $S^3$ with
non-vanishing  modular class.

2. Let $(G,\om)$ be a symplectic Lie group ( for example the
affine group $GA(n)$) ( see [6]). Let $P\times G\too P$ a locally
free action of $G$ on a differentiable manifold $P$ whose
associated foliation will be denoted by $\cal F$. The symplectic
form $\om$ gives arise to a tangential 2-form on $\cal F$ which is
symplectic on the restriction to any leaf of $\cal F$. This gives
canonically  a Poisson structure $\pi$ on $P$ whose symplectic
foliation is $\cal F$.

If the action of $G$ leaves invariant a Riemannian metric on $P$,
the foliation is Riemannian and the modular class of $(P,\pi)$
vanishes. This is the case if $P$ a Lie group and $G$ is a Lie
subgroup which acts by left translations.

Another interesting case is the case where $H$ is a Lie group with
$G$ as subgroup and $\Ga$ is a discrete subgroup of $H$ such that
there is a Riemann metric on $H$ which is right $G$-invariant and
left $\Ga$-invariant. Then the natural homogenous  action of $G$
on $P=\Ga\backslash H$ is locally free and the associated
foliation is Riemannian. For example, the natural action of
$\reel^{2p}$ on the torus $T^n=\ent^n\backslash\reel^n$.

3. Let $G$ be an unimodular Lie group and $H$ an unimodular
subgroup of $G$. Let $\Ga$ a discrete subgroup of $G$. The
homogenous left action of $\Ga$ on $G/H$ lives invariant a volume
form on $G/H$. Let $(M,\om)$ a symplectic manifold with a
 free proper symplectic action $M\times\Ga\too M$. The manifold
$P=M\times_{\Ga}G/H$ carries a foliation $\cal F$ whose leaves are
of the form $M\times\{.\}$. The symplectic form $\om$ gives arise
to a tangential symplectic form and then a Poisson structure
$\pi$. Although the foliation is not Riemannian, the modular class
of $\pi$ vanishes.

4. The affine group $GA(2)$ can be considered as a subgroup of two
simply connected Lie groups od dimension 3 both having a cocompact
discrete subgroup $\Ga$. The first one is $\wi{SL}(2,\reel)$ the
universal covering of $SL(2,\reel)$ and the second is $G_3$ whose
Lie algebra is given by the relations
$$[e_1,e_2]=-e_1,\quad[e_1,e_3]=0,\quad [e_2,e_3]=-e_3.$$
There is Poisson structure ( constructed as above) on the compact
3-manifold $M=\Ga\backslash H$ ( where $H=\wi{SL}(2,\reel)$ or
$G_3$) whose the symplectic foliation is the foliation given by
the homogenous action of $GA(2)$. This foliation ( of codimension
one) is not Riemannian and so the modular class of the Poisson
structure don't vanish.

\section{The first leafwise cohomology space is an invariant of Morita equivalence}

Following [17], recall that a full dual pair
$P_1\stackrel{\rho_1}\leftarrow W\stackrel{\rho_2}\rightarrow P_2$
consists of two Poisson manifolds $(P_1,\pi_1)$ and $(P_2,\pi_2)$,
a symplectic manifold $W$, and two submersions $\rho_1:W\too P_1$
and $\rho_2:W\too P_2$ such that $\rho_1$  is Poisson, $\rho_2$ is
anti-Poisson, and the fibers of $\rho_1$ and $\rho_2$ are
symplectic orthogonal to each other. A Poisson ( or anti-Poisson)
mapping is said to be complete if the pull-back of a complete
Hamiltonian flow under this mapping is complete. A full dual pair
is called complete if both $\rho_1$ and $\rho_2$ are complete. The
Poisson manifolds $P_1$ and $P_2$ are Morita equivalent if there
is exists a complete full dual pair
$P_1\stackrel{\rho_1}\leftarrow W\stackrel{\rho_2}\rightarrow P_2$
such that $\rho_1$ and $\rho_2$ both have connected and simply
connected fibers. Morita equivalent Poisson manifolds $P_1$ and
$P_2$ have isomorphic first Poisson cohomology spaces. More
explicitly, there is a natural isomorphism
$$E:H^1_\pi(P_1)\stackrel{\simeq}\too H^1_\pi(P_2)\eqno(18)$$
which is defined by ( see [4] Lemma 5.2))
$$E([\xi_1])=[\xi_2]\ssi\exists F\in C^\infty(W)/ \mbox{}\quad \xi_1=(\rho_1)_*X_F, \xi_2=-(\rho_2)_*X_F.\eqno(19)$$

Let $\xi_1$ be a representant of a class in $\pi_1^*(H^1_{\cal F}(P_1))$. Then there exists a differential 1-form $\al$ on $P$ such that $\pi(\al)=\xi_1$. Let $\xi_2$ a representant of $E([\xi_1])$. $[\xi_2]\in \pi_2^*(H^1_{\cal F}(P_2))$ if and only if $\xi_2$ is tangent to the symplectic foliation which is true if, for any local Casimir function $f$, $df(\xi_2)=0.$

Let $f$ be a Casimir function. Remark that the relation $(\rho_1)_*X_F=\xi_1$ is equivalent to the existence of $\rho_1$-vertical vector field $V$ such that
$$dF=\rho_1^*(\al)+\om^{-1}(V)$$where $\om:T^*W\too TW$ is the identification associated to the symplectic form on $W$.

We have
\begin{eqnarray*}
df(\xi_2)&=&-df((\rho_2)_*X_F)\\
&=&-d(f\circ\rho_2)(X_F)\\
&=&dF(X_{f\circ\rho_2})\\
&=&\rho_1^*(\al)(X_{f\circ\rho_2})+\om(V,X_{f\circ\rho_2}).
\end{eqnarray*}
$X_{f\circ\rho_2}$ is $\rho_2$-vertical since $f$ is a Casimir function so
$\om(V,X_{f\circ\rho_2})=0$ and $X_{f\circ\rho_2}$ is $\rho_1$-vertical
so $\rho_1^*(\al)(X_{f\circ\rho_2})=0$. So we can conclude.
\bigskip

{\bf References}\bigskip

[1] C. Deninger and   W. Singhof, Real Polarizable Hodge
Structures Arising from Foliations, Annals of Global Analysis and
Geometry 21 (2002) 377-399.

[2] A. El Kacimi, Sur la cohomologie feuillet\'ee, Compositio
Mathematica 49 (1983) 195-215.

[3] V. L. Ginzburg, A. Golubev, Holonomy on Poisson Manifolds and
the modular class, Preprint 1998, math.DG/9812153.

[4] V. L. Ginzburg, J. H. Lu, Poisson cohomology of Morita
equivalent Poisson manifolds, IMRN, 10 (1992), 199-205.

[5] S. Hurder, The Godbillon measure of amenable foliations. J.
Diff. Geom. 23 (1986) 347-365.

[6] A. Medina and P. Revoy, Groupes de Lie à structure
symplectique invariante, in Symplectic Geometry, Groupoids and
integrable systems. MST Publications (1989).

[7] P. Molino, Riemannian foliations, Progress in Mathematics Vol.
73, Birkh\"auser Boston Inc. 1988.

[8] B. Reinhart and J.W. Wood, A metric formula for the
Godbillon-Vey invariant for foliations,
 Proceeddind for the American Society, Vol. 38, Number, April 1973.

[9] P. Tondeur, Foliations on Riemannian manifolds.
Springer-Verlag.

[10] I. Vaisman, Lectures on the Geometry of Poisson Manifolds,
Progress in Mathematics, vol. 118, Birkh\"auser, Berlin, 1994.

[11] A. Weinstein, The Modular Automorphism Group of a Poisson
Manifold, J. Geom. Phys. 23, p. 379-394, 1997.

[12] A. Weinstein, The Local Structure of Poisson Manifolds, J.
Differential Geometry 18, p. 523-557, 1983.

\bigskip
 {\it Abdelhak Abouqateb

  Facult\'e des Sciences et Techniques B.P.
549 Gueliz\\Marrakech. Maroc.\\ e-mail abouqateb@fstg-marrakech.ac.ma \\

 Mohamed Boucetta

  Facult\'e des Sciences et Techniques B.P. 549 Gueliz \\
Marrakech. Maroc.\\ e-mail boucetta@fstg-marrakech.ac.ma}

\end{document}